# Alexander Yurkin

A.M. Prokhorov General Physics Institute of the Russian Academy of Sciences
e-mail: yurkin@fpl.gpi.ru

# RAY TRAJECTORIES, BINOMIAL OF A NEW TYPE, AND THE BINARY SYSTEM

**Abstract**

The paper describes a new algorithm of construction of the nonlinear arithmetic triangle on the basis of numerical simulation and the binary system. It demonstrates that the numbers that fill the nonlinear arithmetic triangle may be binomial coefficients of a new type. An analogy has been drawn with the binomial coefficients calculated with the use of the Pascal triangle. The paper provides a geometrical interpretation of binomials of different types in considering the branching systems of rays.

**Introduction**

Laser generators normally consist of active medium, which enhances light, and a resonator consisting of two plane-parallel mirrors of which one is a mirror with 100% reflection and the other mirror is a semitransparent one [1].

The works [2, 3] offer a mirror of a new type: this is a multi-lobe mirror. Such mirrors consist of numerous semitransparent flat planes inclined towards the laser axis at a small angle $\gamma$ and turned symmetrically around the axis. A ray of light reflects conically from the multi-lobe mirror in the form of a hollow cone, or, to be more precise, in the form of numerous lines situated in the lateral side of the cone. Let us assume for our calculations that the mirrors themselves are thin and that all the rays are inclined at small angles towards the axis, i.e., we use for the calculations the Gaussian [4, 5], or paraxial approximation. This is an approach widely used in calculating optical devices, e.g., telescopes. In a plane (two-dimensional) case, the ray that reflects from our mirror would split in two.

In order to describe the processes in lasers fitted with resonators, wave or geometrical-optics (ray) simulations and equivalent light guide schemes are normally used [1]. Let us note that the light energy emitted by the laser is proportionate to the number of rays.

In the laser fitted with a multi-lobe mirror all the frequencies, excited in active medium, would generate, unlike the laser fitted with a plane-parallel mirror in which generation would take place only at resonance frequencies. Therefore, in the work [6] we named the laser resonator, fitted with a multi-lobe mirror, a quasi-resonator.

The work [7] offered a visual geometrical-optics simulation on the basis of consideration of the binomial distribution in order to describe the propagation of light in the laser fitted with a multi-lobe mirror.

The work [8] described the nonlinear arithmetic pyramid and nonlinear arithmetic triangles and offered recurrence formulas for the construction of the nonlinear arithmetic pyramid and nonlinear arithmetic triangles. The same work shows their physical interpretation, i.e., the description of the process of the propagation of light in lasers fitted with a multi-lobe mirror, and gives examples from the theory of combinations.

In this work we are going to conduct a more detailed research of the nonlinear arithmetic triangle. We will offer a new algorithm of construction of such nonlinear arithmetic triangle and show the relation of this offered algorithm, the binary system and binomial coefficients of different types. The paper provides a geometrical interpretation of binomials of different types in considering ray trajectories for a two-dimensional case.

1. **Pascal triangle and a table algorithm of calculation of binomial coefficients**

1.1. **Pascal triangle**

It is known from [9] that it is possible to calculate binomial coefficients $\binom{n}{p}$

with the use of a two-dimensional numerical table – the Pascal arithmetic triangle in which numbers are positioned by layers in rows (Fig. 1):

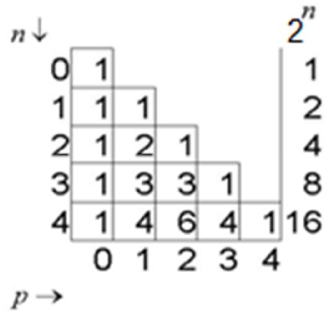

Fig. 1. Pascal arithmetic triangle (linear), $n$– the ordinal number of a row, $p$– the ordinal number of a number in the row, the right column shows the sums of numbers in the row.

The Pascal triangle is constructed in accordance with a recurrence formula:

$$\binom{n}{p} = \binom{n-1}{p-1} + \binom{n-1}{p}, \qquad (1)$$

where $n$ is the binomial's power, $p$ is the ordinal number of a number in the horizontal row, and $0 \le p \le n$.

The numbers in each of the lower rows in the Pascal triangle – binomial coefficients – are the sum of the two numbers of the upper row. The sum of the coefficients for the $n$-power of the binomial equals $2^n$, and is shown in the right column in Fig. 1.

**Example 1.1.** The numbers in the row (Fig. 1) for $n = 3$: 1, 3, 3, 1 are binomial coefficients.

**1.2. Table algorithm of calculation of binomial coefficients**

Let us construct our algorithms on the basis of the following considerations. It is known that in increasing the binomial to power $n$: $(a + b)^n$ we will obtain a polynomial whose elements consist of the multipliers $a$ and $b$ where $a$ and $b$ are certain, for example, real numbers, and $n = 1, 2, ...$.

Let us assume that in a general case addition is communicative while multiplication is non-communicative. Then, for example, in cubing the binomial ($n = 3$) we will obtain:

$$(a + b)^3 =$$
$$aaa_{(0)} + aab_{(1)} + aba_{(2)} + abb_{(3)} + baa_{(4)} + bab_{(5)} + bba_{(6)} + bbb_{(7)}. \quad (2)$$

For descriptive purposes, in (2) we have numbered the elements of the polynomial. The ordinal numbers of $m$ elements of the polynomial, where $m = 0, 1, 2, ...$, are shown in (2) in the form of an inferior index in brackets.

Let us record the obtained polynomial in a tabulated form, with the elements of polynomial (2) positioned as vertical columns in a frame (Fig. 2):

| m → | | | | | | | | k |
|---|---|---|---|---|---|---|---|---|
| 0 | 1 | 2 | 3 | 4 | 5 | 6 | 7 | |
| a | a | a | a | b | b | b | b | 1 ↓ |
| a | a | b | b | a | a | b | b | 2 |
| a | b | a | b | a | b | a | b | 3 |

Fig. 2. Tabulated record of the polynomial (2), with the elements of the polynomial positioned in vertical columns in a frame.

The numbers (Fig. 2) above the table (above the frame) correspond to the ordinal number of $m$-column of the table – an element of polynomial (2), and the numbers on the right of the frame – to the number of $k$-line of the table, where $k = 1, 2, ...$.

Let us replace in Fig. 2 letters $a$ and $b$ with numbers 0 and 1, respectively (Fig. 3):

|   | m → |   |   |   |   |   |   |   |   |   |
|---|---|---|---|---|---|---|---|---|---|---|
|   |   | 0 | 1 | 2 | 3 | 4 | 5 | 6 | 7 |   |
|   | 3 | 0 | 0 | 0 | 0 | 1 | 1 | 1 | 1 | 1 k |
| ↑ | 2 | 0 | 0 | 1 | 1 | 0 | 0 | 1 | 1 | 2 ↓ |
| $l$ | 1 | 0 | 1 | 0 | 1 | 0 | 1 | 0 | 1 | 3 |
|   | Σ | 0 | 1 | 1 | 2 | 1 | 2 | 2 | 3 |   |
| $\binom{n}{p}$ |   | 1 | 3 |   | 3 |   |   |   | 1 |   |

Fig. 3. Tabulated record of numbers in a binary form (binary system), with the numbers positioned in vertical columns in a frame and consisting of 0 and 1.

The sequence of unities and zeros in vertical columns in Fig. 3 may be considered as numbers, recorded in a binary system, and the same numbers $m$ in a denary system are shown in Fig. 3 above the table (above the frame). These numbers also coincide with the numbers above the table (above the frame) in Fig. 2. Numbers $k$ on the right of the frame correspond to a number of the table's line. Numbers $l$ on the left of the frame correspond to number positions in a binary system, where $l = n - k + 1$, in our case number positions $l$ grow bottom-to-top (Fig. 3).

Let us sum up the numbers (i.e., we will sum up the sequence of unities and zeros) in vertical columns (Fig. 3) – their sums $\Sigma$ are shown as bold figures in the first line under the frame. The second line under the frame shows the frequency with which the values of such sums occur: thus, zero occurs once, unity – three times, two – three times, and three – once.

These numbers (in the second line under the frame in Fig. 3), $1, 3, 3, 1$, equal binomial coefficients $\binom{n}{p}$, which, as shown above in Example 1.1, can also be calculated with the use of expression (1) for the Pascal triangle shown in Fig. 1.

If multiplication is ncommunicative, certain (similar) elements of polynomial (2) may be added up in accordance with the values of binomial coefficients $\binom{n}{p}$. In this example ($n = 3$) these are $m = 1, 2, 4$ and $m = 3, 5, 6$, which make up vertical columns $1, 2, 4$ and $3, 5, 6$ of the table in Fig. 2 (and in

accordance with the table in Fig. 3); hence in expression (2) similar elements are:

$$aab_{(1)} = aba_{(2)} = baa_{(4)} = a^2b = A,$$
$$abb_{(3)} = bab_{(5)} = bba_{(6)} = ab^2 = B. \qquad (3)$$

In cubing binomial (2), and with account of (3), we have:

$$(a+b)^3 = a^3 + 3A + 3B + b^3. \qquad (4)$$

Let us name the numbers positioned in horizontal rows of the Pascal triangle (Fig. 1) binomial coefficients of the $p$-type.

Thus, using a tabulated record (Fig. 2, 3) in sec. 1.2, we have obtained a visual demonstration of the elements of the polynomial (Fig. 2) and calculated the binomial coefficients (Fig. 3) without using the image of the Pascal triangle (Fig. 1) or recurrence formula (1). By using the described table algorithm, we have also obtained an expression for binomial (4) from expression (2).

## 2. Non-linear arithmetic triangle and an algorithm of calculation of binomial coefficients of a new type

### 2.1. Non-linear arithmetic triangle

The work [8] described another two-dimensional numerical table – the non-linear arithmetic triangle in which numbers $\begin{bmatrix}n\\q\end{bmatrix}$ are positioned by layers in rows (Fig. 4):

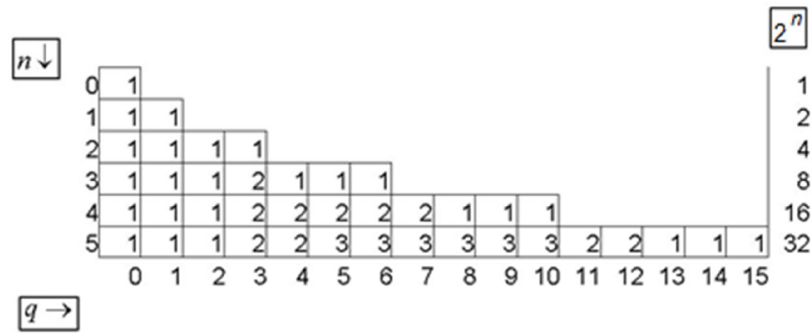

Fig. 4. Non-linear arithmetic triangle, $n$ – the ordinal number of a row, $q$ – the ordinal number of a number in the row, the right column shows the sums of numbers in the row.

The non-linear arithmetic triangle is constructed in accordance with formula [8]:

$$\begin{bmatrix} n \\ q \end{bmatrix} = \begin{bmatrix} n-1 \\ q \end{bmatrix} + \begin{bmatrix} n-1 \\ q-n \end{bmatrix}, \qquad (5)$$

where $n$ is the ordinal number of a row, $q$ is the ordinal number of a number in the horizontal row (Fig.4), and $0 \leq q \leq n(n+1)/2$. The sum of the coefficients for $n$-row, which equals $2^n$, is shown in the right column in Fig. 4.

Expression (5) corresponds to a well-known formula [10] of partitioning numbers into summands. Number $\begin{bmatrix} n \\ q \end{bmatrix}$ shows how many ways there are to partition number $q$ into summands. Each of such summands equals one of numbers $1, 2, 3, \ldots, n$, with no account taken of the position of the summands and all of the summands being different.

**Example 2.1.** Number 3 may be partitioned into the summands $1, 2, 3$ in the two ways: $3 = 1 + 2$, $3 = 3$, i.e., for our case, $n = 3, q = 3$, and there will be: $\begin{bmatrix} 3 \\ 3 \end{bmatrix} = 2$ (Fig. 4).

**2.2. Table algorithm of calculation of binomial coefficients of a new type**

Let us replace in Fig. 2 letters $a$ and $b$ with numbers $0$ and $l$, respectively, i.e., in the first line ($k = 1$) we will have three ($l = 3$) instead of $b$, in the second line ($k = 2$) two ($l = 2$) instead of $b$, and in the third line ($k = 3$) unities ($l = 1$) instead of $b$ (Fig. 5):

|   | $m \rightarrow$ |   |   |   |   |   |   |   |   |
|---|---|---|---|---|---|---|---|---|---|
|   | 0 | 1 | 2 | 3 | 4 | 5 | 6 | 7 |   |
| 3 | 0 | 0 | 0 | 0 | 3 | 3 | 3 | 3 | 1 $k$ |
| $\uparrow$ 2 | 0 | 0 | 2 | 2 | 0 | 0 | 2 | 2 | 2 $\downarrow$ |
| $l$ 1 | 0 | 1 | 0 | 1 | 0 | 1 | 0 | 1 | 3 |
| $\Sigma$ | 0 | 1 | 2 | 3 | 3 | 4 | 5 | 6 |   |
| $\begin{bmatrix} m \\ q \end{bmatrix}$ | 1 | 1 | 1 | 2 |   | 1 | 1 | 1 |   |

Fig. 5. Tabulated record of numbers in a quasi-binary form, with the numbers positioned in vertical columns in a frame and consisting of 0 and numbers $l$.

The numbers above the table (above the frame) in Fig. 5 correspond to the number of $m$-column of the table, and the numbers on the right of the frame correspond to the number of $k$-line of the table. Let us name the sequence of zeros and numbers $l$ a record of numbers in a quasi-binary form. Numbers $l$ are shown on left of the frame (Fig. 5).

Let us sum up the numbers (Fig. 5) in vertical columns – their sums $\Sigma$ are shown as bold figures in the first line below the frame (below the table). The second line below the frame shows the frequency with which the values of such sums occur: thus, zero occurs once, unity – once, two – once, three – twice, four – once, five – once, and six – once.

Numbers in the second line below the frame: $1, 1, 1, 2, 1, 1, 1$, equal numbers $\begin{bmatrix} m \\ q \end{bmatrix}$, positioned in the horizontal row for $n = 3$ of the non-linear arithmetic triangle shown in Fig. 4.

**Example 2.2.** The sum of the numbers positioned in the vertical column

$m = 3$ (Fig. 5) equals: $1 + 2 + 0 = 3$. The sum of the numbers positioned in the vertical column $m = 4$ (Fig. 5) equals: $0 + 0 + 3 = 3$. Thus, the sums of the numbers (the second line below the frame in Fig. 5), which equal three, occur twice, or $\begin{bmatrix}3\\3\end{bmatrix} = 2$, $n = 3$, $q = 3$ (compare with Example 2.1).

Above, in the beginning of sec. 1.2, we assumed that, in a general case, in increasing binomial (2) to power $n$ multiplication would be non-communicative. In the particular case now considered multiplication is non-communicative, too.

However, in the case now considered certain elements of polynomial (2) are similar and may be added up (although multiplication is non-communicative) in accordance with the value of numbers $\begin{bmatrix}n\\q\end{bmatrix}$ by analogy with numbers $\binom{n}{p}$. In this example, these are elements $m = 3, 4$ of polynomial (2), which make up vertical columns $3, 4$ of the table in Fig. 2 (and in accordance with the table in Fig. 5), i.e., in expression (2) in this case (unlike the case described in sec. 1.2, expression (3)) similar elements are:

$$abb_{(3)} = baa_{(4)}, \text{ or, in a concise form: } ab^2 = ba^2 = C. \qquad (6)$$

Below, in sec. 3.2 we will demonstrate geometrical interpretation of similar elements of type (6).

In cubing binomial (2), and with account of (6), we have:

$$(a + b)^3 = a^3 + a^2b + aba + 2C + bab + b^2a + b^3. \qquad (7)$$

Here, to have a concise record, we assumed $aaa_{(0)} = a^3$, $aab_{(1)} = a^2b$, $bba_{(6)} = b^2a$, $bbb_{(7)} = b^3$.

Expression (7), which describes the binomial, differs from expression (4). We will name the numbers, positioned in horizontal lines of the non-linear arithmetic triangle in Fig. 4, binomial coefficients of the $q$-type unlike binomial coefficients of the $p$-type in Fig. 1.

We will also name the binomials, described by expressions (4) and (7), binomials of the $p$- and $q$-types, respectively.

In a similar way, we can also calculate binomial coefficients of the $q$-types for different $n$ values and record a formula similar to (7) for different powers of the $q$-type binomials.

## 3. Branching ray system and geometrical interpretation of binomials

### 3.1. Geometrical interpretation of the $p$-type binomial

Let us consider a symmetric image of the Pascal triangle and a branching ray system in Fig. 6:

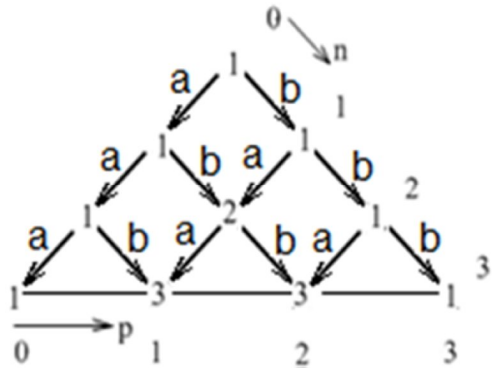

Рис. 6. Symmetric image of the Pascal triangle and a branching ray system, $n$– the number of a ray passage, $p$– the number of branching points (compare with Fig. 1). The ray passage directions are marked with arrows positioned between the numbers of the Pascal triangle: $a$ – arrow direction from the left down and $b$ – arrow direction from the right down.

Let us imagine that the light power (or a particle on the Galton board [11]) propagates (moves) from the top down along the links of broken trajectories. The links (rays) are positioned between the numbers of the Pascal triangle (between the branching points of the rays) in the form of arrows [11] in Fig. 6.

Let us call the links (rays) directed left-down $a$ and those directed right-down – $b$ (Fig. 6).

Thus, for example, for case $n = 3$, there is only one trajectory that leads to point $p = 0$, and this trajectory consists of three links $a$: $aaa$ or $a^3$. Three trajectories, which consist of three links: $aab, aba, baa$, lead to point $p = 1$, and three trajectories lead to point $p = 2$, consisting of links: , etc.

As seen, a portion of light power, which propagates along the rays, comes to the same branching points; therefore, the trajectories by which this power of the rays propagates, are similar ones, and we will obtain from Fig. 6 an expression which is analogous to (3):

$$aab = aba = baa = a^2b = A \text{ and } bba = bab = abb = ab^2 = B.$$

Thus, the light power that propagates along similar trajectories may be added up. Consequently, on the basis of the consideration of the trajectories in Fig. 6, the formula of binomial (2) can be reduced to expression (4).

### 3.2. Geometrical interpretation of the $q$-type binomial

Let us consider another branching system of rays, which originated during the description of processes taking place in lasers. This system was offered by [7] and described in [8] (Fig. 7):

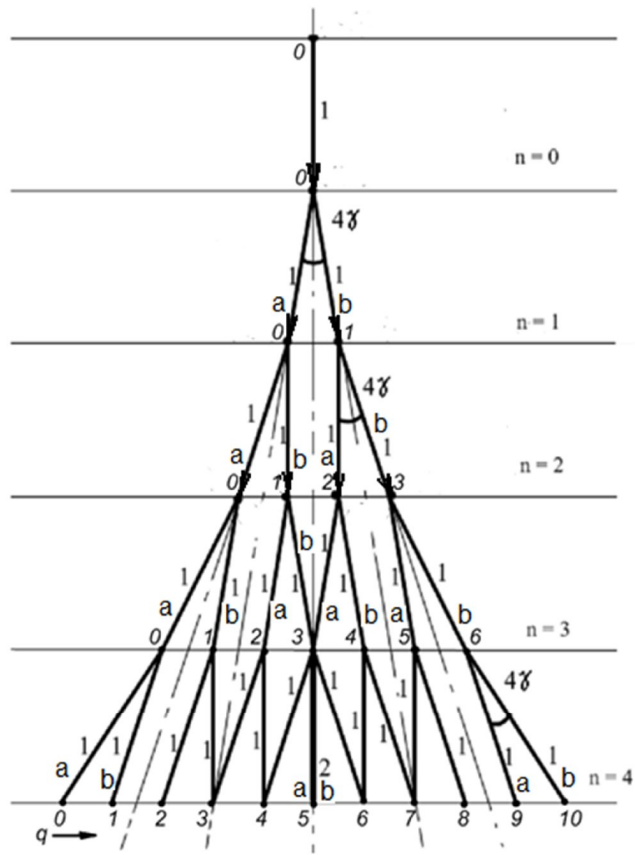

Fig. 7. The branching system of rays inclined at small $2\gamma$ - multiple angles towards the vertical. The rays and links are shown by continuous lines, symmetry axes by line-and-dot lines. Light propagates from the top down, $n$ – number of ray passage, $q$ – number of ray branching points. The rays deflected left are marked $a$, right – $b$.

Let us imagine that light power propagates from the top down in sequence along the links of broken and branching trajectories shown in Fig. 7. (For descriptive purposes, some of the links are shown in the form of an arrow.) The rays that change direction on the left of the direction of the preceding link will be marked $a$, on the right – $b$ (Fig. 7).

Thus, for example, for $n = 3$, there is only one trajectory that leads to point $q = 0$, and this trajectory consists of three links $a$: $aaa$ or $a^3$. There are: one trajectory leading to point $q = 1$, which consists of three links: $aab$; one trajectory leading to point $q = 2$, which consists of three links: $aba$; and two trajectories leading to point $q = 3$: $abb$ and $baa$, etc.

As seen from Fig. 7, a portion of light power, which propagates along the rays, comes to the same branching points; therefore, the trajectories by which this power propagates, are similar ones, and we will obtain an expression which is analogous to (6):

$$abb = baa = ab^2 = ba^2 = C.$$

Thus, the light power that comes to the same branching points may be added up. Consequently, on the basis of the consideration of the trajectories in Fig. 7, the formula of binomial (2) may be recorded as expression (7).

### 3.3. Peculiarities of branching systems

We have considered the two different branching systems of rays in sec. 3.1 and sec. 3.2. It should be noted that the system of rays in Fig. 7 has its own peculiarities as compared to the system of rays in Fig. 6.

The links (rays), making the trajectories in Fig. 6, have one of single-valued directions – $a$ or $b$.

The links (rays), making the trajectories in Fig. 7, may also have one of single-valued directions – $a$ or $b$. This is true, for example, for the rays which make trajectories $aaa = a^3$ or $bbb = b^3$.

The links (rays) in Fig. 7 may also have a non-single-valued direction (in certain cases the same link (ray) has direction $a$ and direction $b$ in other cases), which depends on the direction of the preceding trajectory. This is true, for example, for the rays, which are parallel to the vertical, and for the links which are superimposed (Fig. 7). For example, in the fourth passage ($n = 4$) the ray coming to point $p = 5$ has directions $a$ and $b$ simultaneously, since the last link is a pert of the two trajectories: $abba$ and $baab$ (Fig. 7).

### Conclusion

Sec. 3 demonstrated the geometric interpretation of the $q$-type binomial on the basis of the consideration of trajectories and their links directly in the model shown in Fig. 7. However, the calculation of the trajectory links directly in the figure is not always convenient, although makes it possible to unambiguously determine all the links of which the trajectories consist.

With the use of expression (5) it is possible to construct the non-linear arithmetic triangle (Fig. 4), and, therefore, to calculate how many trajectories in Fig. 7 are similar ones (i.e., come to the same point). However, expression (5) does not show of which particular links the trajectories consist.

We have shown that with the use of the algorithm, offered in sec. 2.2, it is possible not just to calculate the elements (numbers) making up the non-linear arithmetic triangle (i.e., the number of trajectories), but also determine these trajectories themselves and the links of which they consist.

Thus, we have shown the correspondence of the offered algorithm to the geometric model. This is an opportunity provided by the consideration of binomial coefficients and the $q$-type binomial. We have obtained expression (7) for the $q$-type binomial in considering supplementary conditions taken from the processes which take place in laser.

**References**

1. *Yu..A. Ananiev.* Optical resonators and laser beams. Moscow: Nauka, 1990.
2. *A. V. Yurkin.* Patent No. 1777656, USSR, Priority on 26.06.1990; BI No. 43, 1990.
3. *A. V. Yurkin.* New mirror of the laser resonator. // Quantum Electronics, 1991, vol.18, p.393.
4. *V. Bueller.* Gauss. Moscow: Nauka, 1989, p. 149.
5. *G.S. Landsberg.* Optics. Moscow: Nauka, 1976, p. 272.
6. *A. V. Yurkin.* Quasi-resonator – a new interpretation of scattering in laser. // Quantum Electronics, 1994. vol.21, p.385.
7. *A. V. Yurkin.* System of rays in lasers and a new feasibility of light coherence control. // Optics Communications. 1995. V.114, p.393.
8. *A. V. Yurkin.* The ray system in lasers, non-linear arithmetic pyramid and non-linear arithmetic triangles. // Works of Institute of Systems Analysis of Russian Academy of Sciences, 2008, vol. 32(2), pp. 99 – 112.

9. *O.V. Kouzmin.* Pascal triangle: properties and generalizations. SEJ, 2000.
10. *N. Ya. Vilenkin.* Combinatorics. Moscow: Nauka, 1969. (Chapter 4, Formula (24)).
11. *A.N. Kolmogorov, I.G. Zhurbenko, A.V. Prokhorov.* Introduction to the theory of probability. Moscow: Nauka, 1995, p. 20.

# ON BINOMIAL DISTRIBUTION
# OF THE SECOND (NONLINEAR) TYPE
# FOR BIG BINOMIAL POWER

Results of numerical calculations of binomial distribution of the second (nonlinear) type for big power of a binomial are given. Difference of geometrical properties of linear and nonlinear arithmetic triangles and envelopes of binomial distributions of the first and second types is drawn. The empirical formula for half-sums of binomial coefficients of the second (nonlinear) type is offered. Comparison of envelopes of binomial coefficients sums is carried out. It is shown that at big degrees of a binomial a form of envelops of these sums are close.

In work [1] the descriptive geometry optics model on the basis of consideration of binomial distribution for the description of propagation of light in the laser was offered. Calculations were carried out when using paraxial (Gaussian) [2] approximations.

In work [3] the nonlinear arithmetic pyramid and arithmetic triangles of the second (nonlinear) type were described and recurrent formulas for creation of a nonlinear arithmetic pyramid and nonlinear arithmetic triangles were offered. In the same work their physical interpretation – the description of process of distribution of light in the lasers equipped with a multilobe mirror is shown, examples from combination theory are given.

In work [4] the new algorithm of creation of a nonlinear arithmetic triangle was offered, communication of the offered algorithm, binary notation and binomial coefficients of different types is shown. Geometrical interpretation of binomials of two types by consideration of trajectories of rays for a two-dimensional case is given.

Numerical calculations of binomial distribution of the second type are given in the present work for big degrees of a binomial. The analogy to symmetric binomial distribution of the first type [5] is drawn.

**1. Pascal's triangle**
**1.1. About geometrical properties of a triangle of Pascal**

Binomial coefficients can be calculated by means of the two-dimensional numerical table – an arithmetic triangle of Pascal in which numbers are located layer-by-layer in rows (Fig. 1):

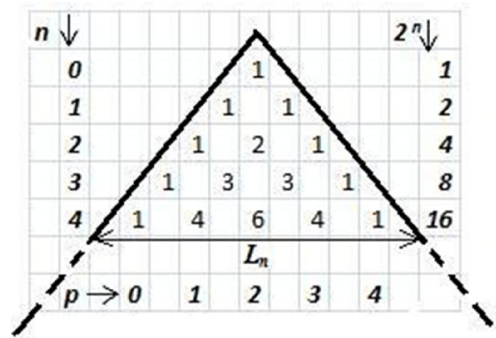

Fig. 1. Pascal's arithmetic triangle (linear), $n$ – serial number of a row, $p$ – serial number of integer in row, in a column on the right – the sums of integers in row, inclined straight lines – side lines of a contour of an arithmetic triangle, $L_n$ – basis length.

Pascal's triangle is under construction according to a recurrent formula:

$$\binom{n}{p} = \binom{n-1}{p-1} + \binom{n-1}{p}, \tag{1}$$

where $n$ – binomial degree, and $p$ – serial number of integer in a horizontal row and

$$0 \leq p \leq n, \quad p_{max} = n.$$

(2)

Numbers in each of the bottom rows in Pascal's triangle – binomial coefficients – are the sum of two numbers of the top row.

We consider that the geometrical sizes of integers (the figures making these integers) are small, and integers in rows are located evenly at identical distances from each other.

The sum of coefficients for $n$ - the binomial degrees, equal $2^n$ is given in a column on the right in Fig. 1. Inclined straight lines in Fig. 1 are elements of a contour of a triangle of Pascal. Length of the basis of an arithmetic triangle $L_n$ for everyone $n$ - a row is proportional to number of binomial coefficient $n = p_{max}$, i.e. $L_n \sim n$. Numbers $n$ - the row, which sum is

$$\sum_{p=0}^{n} \binom{n}{p} = 2^n, \qquad (3)$$

are located within length $L_n$. We consider that binomial coefficients $\binom{n}{p}$ for $n = 4$ (Fig. 1) are located on the line $L_n$ of the basis of a triangle.

Contours of a triangle of Pascal are given in Fig. 2 for big degrees $n$ of a binomial:

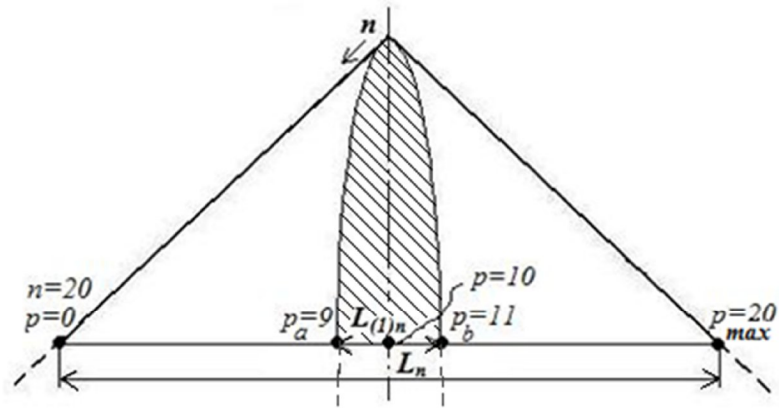

Fig. 2. Contours of a triangle of Pascal for big degrees of a binomial. The numbers filling an arithmetic triangle aren't shown in view of a small scale. For an example points designated binomial coefficients: $\binom{20}{0}$, $\binom{20}{9}$, $\binom{20}{10}$, $\binom{20}{11}$, $\binom{20}{20}$. Inclined straight lines are side lines of a contour of an arithmetic triangle, $L_n$ – length of the basis of a triangle. $L_{(1)\,n}$ – part of the basis of the triangle, concluded between curves.

.

Contours of a triangle of Pascal for great values $n$ in Fig. 2 are constructed similarly by that on Fig. 1. Binomial coefficients for some $n \gg 1$ (Fig. 2) are located in the line of the basis of a triangle.

In $n$ - a row the sum of all binomial coefficients is

$$\sum_{p=0}^{n} \binom{n}{p} = 2^n \sim n \sim L_n. \qquad (4)$$

Part of binomial coefficients $\binom{n}{p_a}, \ldots, \binom{n}{p_b}$ are symmetrized concerning a vertical in the triangle center (the shaded area) in a interval $L_{(1)\,n}$ (Fig. 2).

The sum of this part of binomial coefficients, is approximately equal to a half-sum of (4):

$$\sum_{p_a}^{p_b} \binom{n}{p} = \binom{n}{p_a} + \ldots + \binom{n}{p_b} \approx \frac{1}{2} 2^n$$

(5)

This half-sum $\frac{1}{2} 2^n = 2^{n-1}$ of the binomial coefficients located in the central part in this row, is in limits of length of interval:

$$L_{(1)\,n} \approx L_n^{k_1} \sim \frac{1}{2} 2^n, \qquad (6)$$

and

$$\frac{1}{2} 2^n \sim p_{max}^{k_1} = n^{k_1}, \qquad (7)$$

where (as showed our calculations for $n \gg 1$) empirical coefficient

$$k_1 \approx 0{,}43 \qquad (8)$$

For great values of $n$ it is possible to accept that

$$k_1 \approx \frac{1}{2} \qquad (8a)$$

The interval $L_{(1)\,n}$ in Fig. 2 is located between curves; these curves in a form similar given in [6].

Thus, all binomial (the first type) coefficients (their sum – expression (4)) are located between inclined lines (Fig. 1, 2) and their quantity grows with increase linearly: $L_n \sim n$. Binomial coefficients in the central part Fig. 1, 2 making a half-sum of all binomial coefficients (expression (5, 7)), are located between inclined curve (Fig. 2) and their number grows with increase $n$ according to dependence:

$$L_{(1)n} \approx L_n^{k_1} \approx L_n^{1/2} \sim 2^{n-1} \sim an^{1/2}, \qquad (6a)$$

where $a$ is some positive coefficient.

### 1.2. Binomial distribution of the first type

For convenience of a statement we will use the standard concept "probability". The binomial probability (the relation of some binomial coefficient to the sum of binomial coefficients) is equal $P_n(p) = \binom{n}{p}\frac{1}{2^n}$ [5].

Binomial distribution of the first type [5, 8] and forms of envelopes for binomial distribution are shown in Fig. 3.

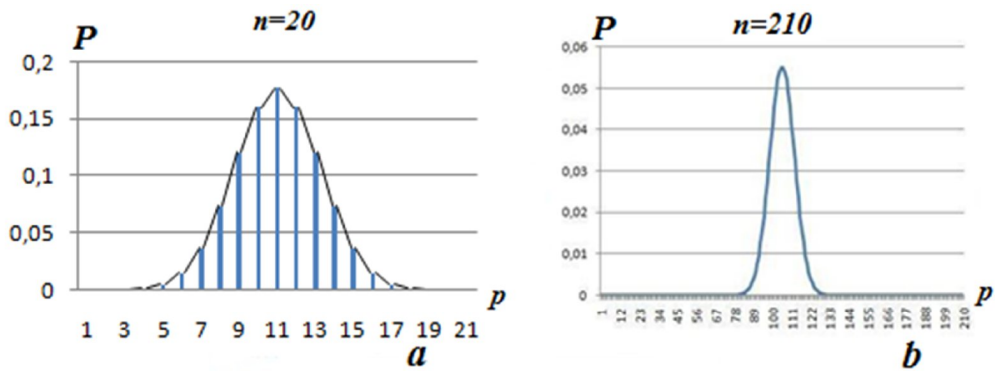

Fig. 3. Binomial distribution of the first type for $n = 20$ (a) and envelopes of binomial

distribution of the first type $P_n(p)$ for $n = 20$ (a) and $n = 210$ (b).

In Fig. 4 the form of envelope (similar to a form for normal distribution [7, 8]) for the sums of binomial coefficients $\sum_{p=0}^{n} \binom{n}{p}$ of the first type is shown.

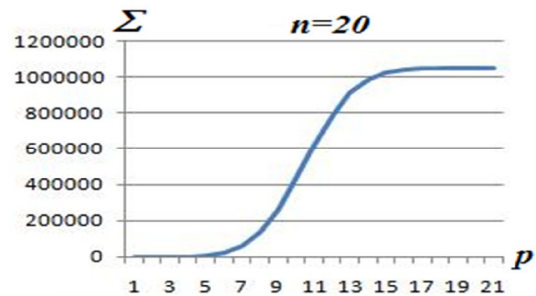

Fig. 4. The envelope of the sums of binomial coefficients of the first type for $n = 20$.

In more detail about binomial distribution of the first type for big pover (Moivre-Laplace's theorem) it is possible to find, for example in [5].

## 2. Nonlinear arithmetic triangle
### 2.1. On geometrical properties of a nonlinear arithmetic triangle

The binomial coefficients of the second type $\begin{bmatrix} n \\ q \end{bmatrix}$ described in works [3, 4], it is possible to calculate by means of other two-dimensional numerical table – a nonlinear arithmetic triangle in which numbers $\begin{bmatrix} n \\ q \end{bmatrix}$ are located layer-by-layer in rows (Fig. 5).

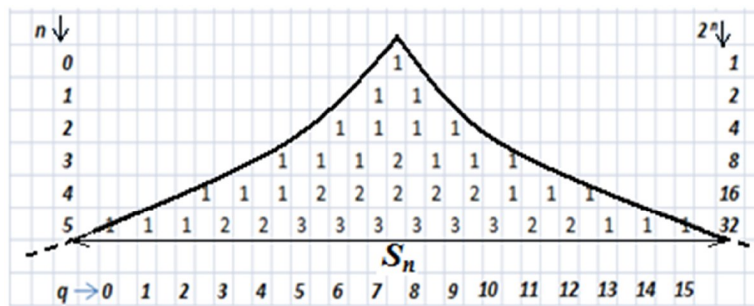

Fig. 5. Nonlinear arithmetic triangle. $n$ – serial number of a row, $q$ – serial number of integer in row, in a column on the right – the sums of integers in row, inclined curves – sides lines of a contour of an arithmetic triangle, $S_n$ – basis length.

The nonlinear arithmetic triangle is under construction according to a formula [3]:

$$\begin{bmatrix} n \\ q \end{bmatrix} = \begin{bmatrix} n-1 \\ q-n \end{bmatrix} + \begin{bmatrix} n-1 \\ q \end{bmatrix}, \tag{9}$$

where $n$ – binomial degree, $q$ – serial number of integer in a horizontal row (Fig. 5) and

$$0 \leq q \leq n(n+1)/2,\ q_{max} = n(n+1)/2. \tag{10}$$

The sum of coefficients for $n$ - the binomial power is given in a column on the right in Fig. 5.

Expression (9) corresponds to a known formula [9] of integers partition. The number $\begin{bmatrix} n \\ q \end{bmatrix}$ shows by how many method the number $q$ may be partitioned into items each of which is equal to one of the numbers $1, 2, 3, \ldots, n$, the order of items being not taken into account and all items being different.

Inclined lines in Fig. 5 are the elements of a contour of a nonlinear arithmetic triangle. Length of the basis $S_n$ of an arithmetic triangle for everyone $n$ - row is proportional $q_{max}$, i.e.

$$S_n \sim q_{max} = n(n+1)/2. \tag{11}$$

Numbers $n$ - row, which sum

$$\sum_{q=0}^{q_{max}} \begin{bmatrix} n \\ q \end{bmatrix} = 2^n \tag{12}$$

are located within length $S_n$. We consider that binomial coefficients $\begin{bmatrix} n \\ p \end{bmatrix}$ for $n = 4$ (Fig. 5) are located on the line $S_n$ of the basis of a triangle.

Contours of a nonlinear arithmetic triangle are given in Fig. 6 for big power $n$ of a binomial.

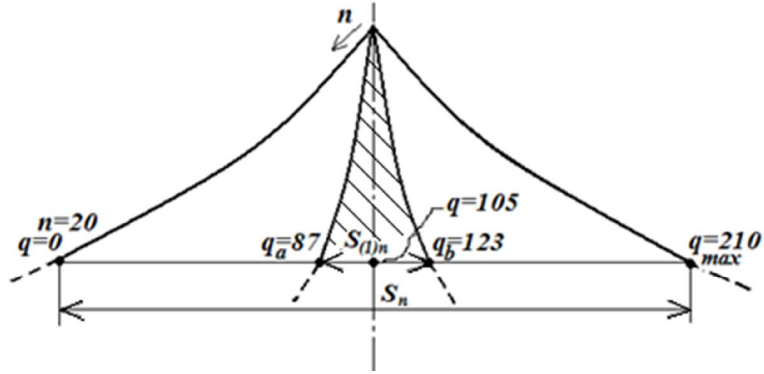

Fig. 6. Contours of a nonlinear arithmetic triangle at big power $n$ of a binomial. The numbers filling an arithmetic triangle, aren't shown in view of a small scale. For an example points designated binomial coefficients: $\begin{bmatrix}20\\0\end{bmatrix}$, $\begin{bmatrix}20\\87\end{bmatrix}$, $\begin{bmatrix}20\\105\end{bmatrix}$, $\begin{bmatrix}20\\123\end{bmatrix}$, $\begin{bmatrix}20\\210\end{bmatrix}$. External inclined curves are side lines of a contour of a nonlinear arithmetic triangle, $S_n$ – length of the basis of a triangle. $S_{(1)\,n}$ – part of the basis of the triangle, concluded between internal (located closer to a vertical) curves.

Contours of a nonlinear arithmetic triangle for great values in Fig. 6 are constructed similarly by that on Fig. 5. Binomial coefficients for some $n \gg 1$ (Fig. 6) are located along on the line $S_n$ of the basis of a triangle.

In $n$ - a row the sum of all binomial coefficients of the second type

$$\sum_{q=0}^{q_{max}} \begin{bmatrix}n\\q\end{bmatrix} = 2^n \sim q_{max} = \frac{n(n+1)}{2} \sim S_n. \qquad (13)$$

Part of binomial coefficients $\begin{bmatrix}n\\q_a\end{bmatrix}, \dots, \begin{bmatrix}n\\q_b\end{bmatrix}$ is symmetrized concerning a vertical in the triangle center on interval $S_{(1)\,n}$ (Fig. 6).

The sum of this part of binomial coefficients is approximately equal to a half-sum from (12):

$$\sum_{q_a}^{q_b} \begin{bmatrix}n\\q\end{bmatrix} = \begin{bmatrix}n\\q_a\end{bmatrix} + \dots + \begin{bmatrix}n\\q_b\end{bmatrix} \approx \frac{1}{2}2^n.$$

(14)

This half-sum $\frac{1}{2}2^n = 2^{n-1}$ of the binomial coefficients located in the central part of a triangle (the shaded area) in the row as showed our numerical calculations, is in interval of length

$$S_{(1)\,n} \approx S_n^{k_2} \sim \frac{1}{2}2^n, \qquad (15)$$

and

$$\frac{1}{2}2^n \sim q_{max}^{k_2} = \left[\frac{n(n+1)}{2}\right]^{k_2}, \qquad (16)$$

where (as showed our calculations for $n \gg 1$) empirical coefficient

$$k_2 \approx 0{,}71 \qquad (17)$$

For great values of $n$ it is possible to accept that

$$k_2 \approx \frac{1}{2} \qquad (17a)$$

The interval $S_{(1)\,n}$ is located in the central part between curves in Fig. 6.

Thus, all binomial coefficients of the second type (their sum - expression (12)) are located between inclined curves of a contour of an arithmetic triangle (Fig. 5, 6) and their quantity grows with increase according to square dependence. Binomial coefficients in the central part of Fig. 6, (the shaded area) components a half-sum (expressions (14, 16, 17, 17a)) are located between inclined curves (Fig. 6) and their number grows with increase according to empirical dependence:

$$S_{(1)n} \approx S_n^{k_2} \approx S_n^{1/2} \sim 2^{n-1} \sim bn \qquad (16a)$$

where $b$ is some positive coefficient.

It should be noted that value empirical coefficients $k_2$ and $b$ demands further specification and possible analytical researches.

### 2.2. Binomial distribution of the second type

Forms of envelopes for binomial distribution of the second type are shown in Fig. 7.

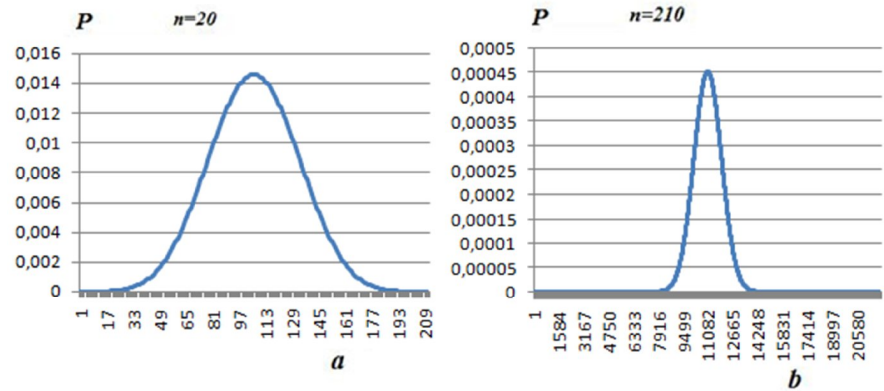

Fig. 7. Envelopes of distribution of the second type $P_{(1)n}(q)$ for $n = 20$ (a) and $n = 210$ (b).

In Fig. 8 the envelope (similar to Fig. 4) for the sums of binomial coefficients of the second type $\sum_{q=0}^{q_{max}} \begin{bmatrix} n \\ q \end{bmatrix}$ is shown.

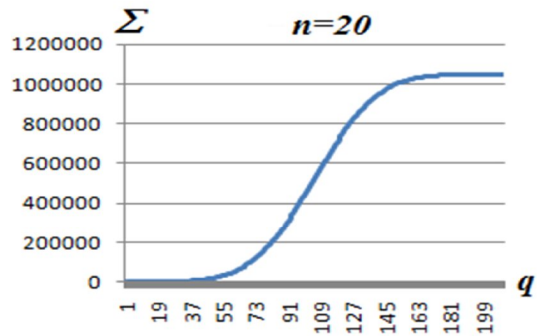

Fig. 8. Envelope of the sums of binomial coefficients of distribution of the second type for $n = 20$.

### 3. Comparison of forms of envelopes of binomial distributions of the first and second types

In Fig. 9 are provided in identical scale envelopes for binomial distributions of first (Fig. 3a) and second (Fig. 7a) of types for $n = 20$.

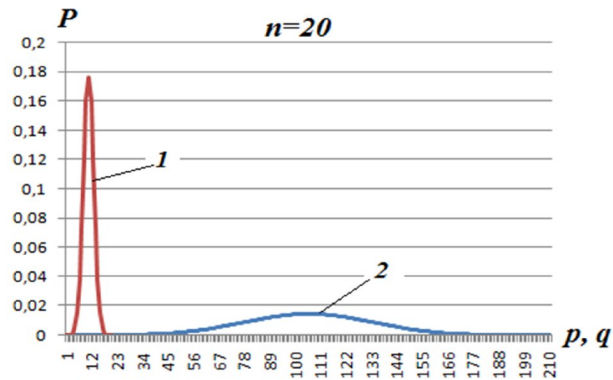

Fig. 9. Envelopes of binomial distributions of the first (1) and the second (2) types for $n = 20$.

In Fig. 10 are provided in identical scale envelopes of the sums of binomial coefficients of first (Fig. 4) and the second (Fig. 8) types for $n = 20$.

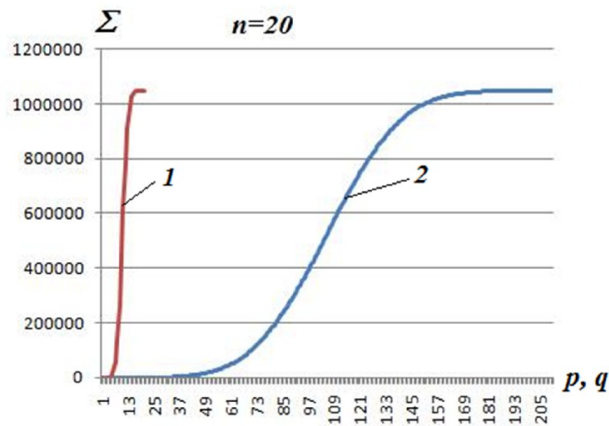

Fig. 10. Envelopes of the sums of binomial coefficients of the first (1) and the second (2) types for $n = 20$.

On Fig. 11 are provided in various (on abscissa axis) scale envelopes of the sums of binomial coefficients of the first and second types for $n = 20$.

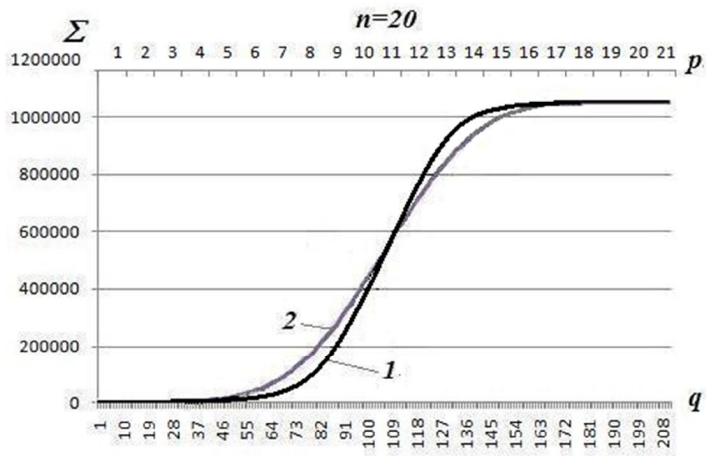

Fig. 11. Envelopes of the sums of binomial coefficients of the first (1) and the second (2) types for $n = 20$. Scale on abscissa axis is various for $p$ and $q$, differs by 10,5 times.

Scale on abscissa axis in Fig. 11 differs by 10,5 times, since $\frac{q_{max}}{p_{max}} = 10,5$ for a case $n = 20$. In Fig. 11 it is visible that the size (3) sums $\sum_{p=0}^{n} \binom{n}{p}$ of binomial coefficients of the first type increases quicker, than the size (12) sums $\sum_{q=0}^{q_{max}} \left[\begin{array}{c} n \\ q \end{array}\right]$ of binomial coefficients of the second type.

As show our numerical calculations, at great values $n$ a form of the sums of binomial coefficients envelopes (Fig. 12) practically don't differ.

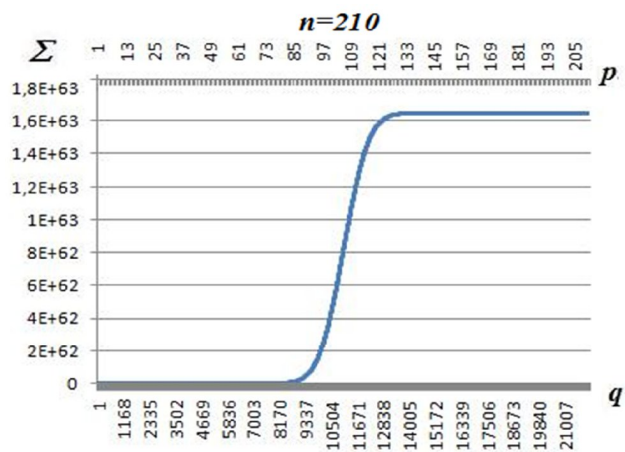

Fig. 12. Envelopes of the sums of binomial coefficients of the first and second types for $n = 210$ (are imposed at each other). Scale on abscissa axis is various for $p$ and $q$, differs by 105,5 times.

Scale on abscissa axis in Fig. 12 differs by 105,5 times, since $\frac{q_{max}}{p_{max}} = 105,5$ for a case $n = 210$. In Fig. 12 it is visible that forms of curves coincide.

It is possible to assume that at great values $n$ ($n > 100$), will be close in a form not only envelopes of the sums (Fig. 12) of binomial coefficients, but also envelopes of binomial distributions of the first and second types (Fig. 3, 7).

**Conclusion**

In the present work numerical calculations and descriptive geometrical interpretations of properties of linear and nonlinear arithmetic triangles were presented.

Compliance of the sums and half-sums of binomial coefficients to certain lengths of intervals of linear and nonlinear triangles was shown.

The empirical formula (demanding further specification) for calculation of half-sums of the binomial coefficients located in the central part of a nonlinear arithmetic triangle is offered.

Forms of envelopes of binomial distributions of the first and second types, and also forms of envelopes for the sums of binomial coefficients were shown.

Comparison of forms of envelopes for sums of binomial coefficients showed that at big degrees of binomials a form of the sums envelopes these (at the certain scales of the image) practically don't differ.

**References**

1. *A. V. Yurkin*. System of rays in lasers and a new feasibility of light coherence control. // Optics Communications. 1995. V.114, p.393.
2. *G. S. Landsberg*. Optics. Moscow: Nauka, 1976, p. 272. (Russian)

3. *A. V. Yurkin*. The ray system in lasers, non-linear arithmetic pyramid and non-linear arithmetic triangles. // Proceedings of Institute of Systems Analysis of RAS, 2008, vol. 32, no. 2, pp. 99 – 112. (Russian) arXiv:1302.5214
4. *A. V. Yurkin*. Ray trajectories and the algorithm to calculate the binomial coefficients of a new type // Proceedings of Institute of Systems Analysis of RAS, 2009, vol. 42, no.1, pp. 66 - 77. (Russian) arXiv:1302.4842

5. *A. N. Kolmogorov, I .G. Zhurbenko, A.V. Prokhorov.* Introduction to the theory of probability. Moscow: Nauka, 1995. (Russian)

6. The mathematical encyclopedic dictionary / under the editorship of *Y. V. Prokhorov*, M: The Soviet encyclopedia, 1988, the Art. of Bernoulli wandering, p. 86. (Russian)

7. *W. Feller.* An introduction to probably theory and its applications. Moscow: Mir, 1964, vol. 1, ch. 7, p. 182. (Russian) [N. Y., John Wiley & Sons Inc.]

8. *G. Korn, T. Korn.* Mathematical handbook for scientists and engineering. Moscow: Nauka, 1974, ch. 18, pp. 572, 576. (Russian) [ N. Y. , 1968]

9. *N. Ia. Vilenkin*. Combinatorics. Moscow: Nauka, 1969. (Ch. 4, formula (24)). (Russian)

**Thanks**